\magnification 1250
\mathcode`A="7041 \mathcode`B="7042 \mathcode`C="7043
\mathcode`D="7044 \mathcode`E="7045 \mathcode`F="7046
\mathcode`G="7047 \mathcode`H="7048 \mathcode`I="7049
\mathcode`J="704A \mathcode`K="704B \mathcode`L="704C
\mathcode`M="704D \mathcode`N="704E \mathcode`O="704F
\mathcode`P="7050 \mathcode`Q="7051 \mathcode`R="7052
\mathcode`S="7053 \mathcode`T="7054 \mathcode`U="7055
\mathcode`V="7056 \mathcode`W="7057 \mathcode`X="7058
\mathcode`Y="7059 \mathcode`Z="705A
\def\spacedmath#1{\def\packedmath##1${\bgroup\mathsurround
=0pt##1\egroup$} \mathsurround#1
\everymath={\packedmath}\everydisplay={\mathsurround=0pt}} 
\def\nospacedmath{\mathsurround=0pt
\everymath={}\everydisplay={} } \spacedmath{2pt}
\font\eightrm=cmr8     \font\sixrm=cmr6   
\font\cc=cmcsc10
\def\pc#1{\tenrm#1\sevenrm}
\def\tx{\kern-1.5pt -}
\def\cqfd{\kern 2truemm\unskip\penalty 500\vrule height
4pt depth 0pt width 4pt} 
\def\ind{\par\hskip 0,8truecm\relax}

\parindent=0cm
\vsize = 25truecm
\hsize = 16truecm
\voffset = -.5truecm
\baselineskip15pt
\overfullrule=0pt
\centerline{\bf The Szpiro inequality for higher genus fibrations}
\smallskip
\smallskip \centerline{Arnaud {\pc BEAUVILLE}}
\vskip1.2cm 

{\bf Introduction}
\smallskip
\ind The aim of this note is to prove the following result:

{\bf Proposition}$.-$ {\it Let $f:S\rightarrow B$ be a
non-trivial semi-stable fibration of genus $g\ge 2$, $N$
the number of critical points of $f$ and $s$ the number of
singular fibres. Then}
$$N< (4g+2)(s+2g(B)-2)\ .$$
\ind Recall that a semi-stable fibration of genus $g$ is a
surjective holomorphic map of a smooth projective surface
$S$ onto a smooth curve $B$, whose generic fibre is a
smooth curve of genus $g$ and whose singular fibres are
allowed only ordinary double points; moreover we impose
that each smooth rational curve contained in a fibre meets
the rest of the fibre in at least 2 points (otherwise by
blowing up  non-critical points of $f$ in a singular fibre 
we could arbitrarily increase $N$ keeping
$s$ fixed).
\ind The corresponding inequality $N\le 6(s+2g(B)-2)$ in
the case $g=1$ has been observed by Szpiro; it was
motivated by the case of curves over a number field, where
an  analogous inequality would have far-reaching
consequences [S]. The higher genus case is considered in
the recent preprint [BKP], where the authors prove the
slightly weaker inequality
$N\le (4g+2)s$ for hyperelliptic fibrations over ${\bf
P}^1$. Their method is topological, and in fact the result
applies in the much wider context of symplectic Lefschetz
fibrations. We will show that in the more restricted
algebraic-geometric set-up, the Proposition is a direct
consequence of two classical inequalities in   surface
theory. It would be interesting to know whether the proof
of [BKP] can be extended to non-hyperelliptic fibrations. 
\vskip1truecm
{\bf Proof}
\ind The main numerical invariants of a surface $S$ are the
square of the canonical bundle $K_S$, the Euler-Poincar\'e
characteristic $\chi({\cal O}_S)$ and the topological
Euler-Poincar\'e characteristic
$e(S)$; they are linked by the Noether formula
$12\,\chi({\cal O}_S)=K_S^2+e(S)$. For a semi-stable
fibration
$f:S\rightarrow B$ it has become customary to modify these
invariants as follows. Let $b$ be the genus of $B$, and
$K_{f}=K_X\otimes f^*K_B^{-1}
$  the relative canonical bundle of $X$ over $B$; then we
consider:
$$\nospacedmath\displaylines{K_f^2=K_X^2-8(b-1)(g-1)\cr
\chi_f:=\deg f_*(K_f)=\chi({\cal O}_X)-(b-1)(g-1)\cr
e_f:=N=e(X)-4(b-1)(g-1)\ .}$$
Observe that we have again $12\,\chi_f=K_f^2+e_f$. We will
use the Xiao inequality ([X], Theorem 2)
$$K_f^2\ge (4-{4\over g})\,\chi_f$$
and the ``strict canonical class inequality" ([T], lemma
3.1)
$$ K_f^2< 2(g-1)(s+2b-2)\quad {\rm for}\  \ s>0\ .
$$ 

\ind Let us prove the Proposition. If $s=0$, we have
$N=0$ and $g(B)\ge 2$ (otherwise the fibration would
be trivial), so the inequality of the Proposition holds.
Assume
$s>0$; the Xiao inequality gives
$$3g\,K_f^2\ge 12(g-1)\,\chi_f=(g-1)(K_f^2+e_f)\ ,$$
hence, using the strict canonical class inequality,
$$N=e_f\le\  {2g+1\over g-1}\ K_f^2 < 
(4g+2)(s+2b-2)\ .\cqfd$$
\bigskip
{\it Example}$.-$ We constructed in [B]  a
semi-stable genus 3 fibration over ${\bf P}^1$ with 5
singular fibres; each of these has $8$ double points.
Therefore
$$N=40\quad{\rm  and}\quad (4g+2)(s-2)=42\ .$$ 
\vskip1.5cm
\centerline{REFERENCES} \vglue15pt\baselineskip12.8pt
\def\num#1{\smallskip\item{\hbox to\parindent{\enskip 
[#1]\hfill}}}\parindent=1.3cm
\num{B} {\cc A. Beauville}: {\sl 	Le nombre minimum de
fibres singuli\`eres d'une courbe stable sur ${\bf P}^1$.} 
Ast\'erisque {\bf 86}  (1981), 97-108 . 
\num{BKP} {\cc F. Bogomolov, L. Katzarkov, T. Pantev}:
 {\sl Hyperelliptic Szpiro inequa\-lity}.
Preprint math.GT/0106212.
\num{S}  {\cc L. Szpiro}: {\sl Discriminant et conducteur
des courbes elliptiques}.  S\'eminaire sur les Pinceaux
de Courbes Elliptiques (Paris, 1988). Ast\'erisque {\bf
183} (1990), 7--18. 
\num{T} {\cc S.-L.  Tan}: {\sl The minimal number of
singular fibers of a semistable curve over} ${\bf P}^1$.
J. Algebraic Geom. {\bf 4} (1995),  591--596. 
\num{X} {\cc G. Xiao}: {\sl Fibered algebraic surfaces with
low slope}. Math. Ann. {\bf 276} (1987),  449--466. 
\vskip0.8cm
\def\pc#1{\eightrm#1\sixrm}
\hfill\vtop{\baselineskip12pt\eightrm\hbox to 5cm{\hfill
Arnaud {\pc BEAUVILLE}\hfill}
 \hbox to 5cm{\hfill Laboratoire J.-A. Dieudonn\'e\hfill}
 \hbox to 5cm{\sixrm\hfill UMR 6621 du CNRS\hfill}
\hbox to 5cm{\hfill {\pc UNIVERSIT\'E DE}  {\pc NICE}\hfill}
\hbox to 5cm{\hfill  Parc Valrose\hfill}
\hbox to 5cm{\hfill F-06108 {\pc NICE} Cedex 02\hfill}}
\end